\documentclass[12pt,a4paper]{article}
\usepackage[latin1]{inputenc}
\usepackage[T1]{fontenc}
\usepackage{amsmath,graphics,amssymb,amsthm,amsfonts,latexsym,color,mathtools}

\usepackage{natbib}

\newtheorem{theorem}{Theorem}

\DeclareMathOperator*{\argmin}{arg\,min}

\begin{document}
\begin{center}
{\large {\bf On oracle efficiency of the ROAD classification rule} 
\vspace{3mm} \\    
Britta Anker Bak and 
Jens Ledet Jensen}
\vspace{2mm} \\
Department of Mathematics, Aarhus University, Denmark
\end{center}

\author{BRITTA ANKER BAK}
\author{JENS LEDET JENSEN}

\section*{Abstract}

For high-dimensional classification Fishers rule performs poorly due to 
noise from estimation of the covariance matrix. \cite{Fan} 
introduced the ROAD classifier that puts an $L_1$-constraint on the 
classification vector. In their Theorem 1 \cite{Fan} 
show that the ROAD classifier asymptotically has the same 
misclassification rate as the corresponding oracle based classifier. 
Unfortunately, the proof contains an error. Here we restate the theorem 
and provide a new proof. 

\section{Introduction}

We consider classification among two groups based on a $p$-dimensional 
normally distributed variable. Let the means in the two groups 
be $\mu_1$ and $\mu_2$, 
and let the common variance be $\Sigma$. Also, let 
the probability of belonging to either of the two groups be $\frac{1}{2}$. 
Defining $\mu_a=(\mu_1+\mu_2)/2$ and $\mu_d=(\mu_1-\mu_2)/2$, 
the Bayes discriminant rule becomes
\[
 \delta_w(x)=1+1(w^T(x-\mu_a)<0),\ \ \text{with}\ w=w_F=\Sigma^{-1}\mu_d,
\]
where $x$ is classified to group 1 or 2 according to the value of 
$\delta_w(x)$. 
The misclassification rate of the rule $\delta_w$ is 
\[
 W(\delta_w)=\bar \Phi\bigl(\frac{1}{2}w^T\mu_d/(w^T\Sigma w)^{1/2}\bigr),
\]
where $\bar \Phi(z)=1-\Phi(z)$ is the upper tail probability of a 
standard normal distribution. 
The interpretation of the Bayes rule is that $w_F$ is the vector 
that minimizes the misclassification rate. \cite{Fan} suggest 
to use a $L_1$ regularized version of $w_F$, that is, 
\[
 w_c=\argmin_{\|w\|_1\leq c,\ 
 w^T\mu_d=1} w^T\Sigma w.
\]
Its sample version
\[
 \hat w_c=\argmin_{\|w\|_1\leq c,\ 
  w^T\hat\mu_d=1} w^T\hat\Sigma w
\]
yields the ROAD classifier
\[
 \hat\delta=1+1(\hat w_c^T(x-\hat\mu_a)<0).
\]

Theorem 1 of \cite{Fan} states that the misclassification rate 
$W(\hat\delta)$ of the ROAD classifier approaches the 
misclassification rate of the 
oracle classifier $W(\delta_{w_c})$. Unfortunately, an essential 
step in the proof use an inequality which is not valid, 
see Appendix~A for details. We reformulate the theorem and give 
a new proof. 

\begin{theorem}
Let $\epsilon$ be a positive constant such that 
$\max_j\{|\mu_{dj}|\}> \epsilon$, and 
$c>\epsilon+1/\max_j\{|\mu_{dj}|\}$. Let $a_n$ be a 
sequence tending to zero such that 
$\| \hat\Sigma-\Sigma \|_\infty=O_p(a_n)$, and 
$\| \hat\mu_i-\mu_i \|_\infty=O_p(a_n)$, $i=1,2$. 
Then, as $n\rightarrow\infty$:
\[
 W(\hat\delta)-W(\delta_{w_c})=O_p(d_n)
\]
with $d_n=c^2 a_n(1+c^2\|\Sigma\|_\infty)$.
\end{theorem}

Prior to proving the theorem we comment on the differences compared to 
Theorem 1 of \cite{Fan}. Contrary to us, \cite{Fan} requires that the smallest 
eigenvalue of $\Sigma$ is bounded from below. 
The upper bound on $W(\hat\delta)-W(\delta_{w_c})$ 
in \cite{Fan} depends on the sparsity of $w_c$ and 
of $w_c^{(1)}$, where $w_c^{(1)}$ is given by
\[
  w_c^{(1)}=\argmin_{\|w\|_1\leq c,\ 
  w^T\hat\mu_d=1} w^T\Sigma w,
\]
whereas our bound depends on the regularizing parameter $c$ only. 
In the formulation of the theorem $c$ is allowed to depend on $n$. 
We require a lower bound on $\max_j\{|\mu_{dj}|\}$, which is not 
part of the theorem in \cite{Fan}. However, it enters indirectly 
in that we must have $c>1/\max_j\{|\mu_{dj}|\}$ in order for $w_c$ to exist. 
Thus, if $\max_j\{|\mu_{dj}|\}\rightarrow 0$, we have $c\rightarrow\infty$, 
and $c$ enters the upper bound of \cite{Fan}. The reason for our more restrictive condition $c>\epsilon+1/\max_j\{|\mu_{dj}|\}$ 
is that the theorem only makes sense 
if $\hat w_c$ exists with probability tending to one. 
Similarly, whereas \cite{Fan} have the condition 
$\| \hat\mu_d-\mu_d \|_\infty=O_p(a_n)$,
we have
$\| \hat\mu_i-\mu_i \|_\infty=O_p(a_n)$, $i=1,2$, 
in order to handle a term in 
the misclassification rate that has been neglected in \cite{Fan}. 
Finally, $\|\Sigma\|_\infty$ appears in our bound. However, requiring 
that the variances $\Sigma_{ii}$, $i=1,\ldots,p$, are bounded is 
often encountered in high dimensional settings. 


\section{Proof of Theorem 1}

In the proof we use the following inequalities:
\begin{align}
 | \bar\Phi(a(1+\epsilon) -\bar\Phi(a)|
 & \leq 2\epsilon\ \text{for}\ a>0\ \text{and}\ |\epsilon|<1,
 \label{1}\\
 | \bar\Phi((a+\epsilon)^{-1/2})-\bar\Phi(a^{-1/2}) | 
 & \leq \epsilon \ \text{for}\ a>0\ \text{and}\ a+\epsilon>0.
 \label{2}
\end{align}
The misclassification rate consists of two terms corresponding 
to an observation from each of the two groups. 
The proofs for the two terms are identical, 
so to simplify we consider the misclassification rate of an observation from group 1 only. 
Using \eqref{1} the misclassification rate of $\hat\delta$ becomes 
\begin{align}
 W(\hat\delta) &=
 \bar\Phi\Bigl( \frac{1}{2}
 \frac{\hat w_c^T\hat\mu_d+\hat w_c^T(\hat\mu_1-\mu_1)}
 {\sqrt{\hat w_c^T\Sigma\hat w_c}}
 \Bigr)
 = 
 \bar\Phi\Bigl( \frac{1}{2}
 \frac{1}{\sqrt{\hat w_c^T\Sigma\hat w_c}}
 \Bigr) + O(|\hat w_c^T(\hat\mu_1-\mu_1)|)
 \nonumber
 \\
 & \leq \bar\Phi\Bigl( \frac{1}{2}
 \frac{1}{\sqrt{\hat w_c^T\Sigma\hat w_c}}
 \Bigr) + O(c\| \hat\mu_1-\mu_1\|_\infty).
 \label{3}
\end{align}
Next,
\[
 | \hat w_c^T\Sigma\hat w_c - \hat w_c^T\hat\Sigma\hat w_c  |
 \leq c^2\| \hat\Sigma-\Sigma \|_\infty,
\]
and from \eqref{2} we get
\begin{equation}
  \bar\Phi\Bigl( \frac{1}{2}
 \frac{1}{\sqrt{\hat w_c^T\Sigma\hat w_c}}
 \Bigr) = 
  \bar\Phi\Bigl( \frac{1}{2}
 \frac{1}{\sqrt{\hat w_c^T\hat\Sigma\hat w_c}}
 \Bigr)
 + O(c^2\| \hat\Sigma-\Sigma \|_\infty).
 \label{4}
\end{equation}
From the proof in \cite{Fan} we see that
\[
 | \hat w_c^T\hat\Sigma\hat w_c - w_c^{(1)T}\Sigma w_c^{(1)}  |
 \leq c^2\| \hat\Sigma-\Sigma \|_\infty,
\]
and thus
\begin{equation}
  \bar\Phi\Bigl( \frac{1}{2}
 \frac{1}{\sqrt{\hat w_c^T\hat \Sigma\hat w_c}}
 \Bigr) = 
  \bar\Phi\Bigl( \frac{1}{2}
 \frac{1}{\sqrt{ w_c^{(1)T}\Sigma\hat w_c^{(1)}}}
 \Bigr)
 + O(c^2\| \hat\Sigma-\Sigma \|_\infty).
 \label{5}
\end{equation}
Combining (\ref{3}--\ref{5}) we have
\begin{equation}
 W(\hat\delta) = 
  \bar\Phi\Bigl( \frac{1}{2}
 \frac{1}{\sqrt{ w_c^{(1)T}\Sigma\hat w_c^{(1)}}}
 \Bigr)
 + O(c^2\| \hat\Sigma-\Sigma \|_\infty + 
     c\| \hat\mu_1-\mu_1\|_\infty).
 \label{6}
\end{equation}
Since the oracle misclassification rate is 
$W(\delta_{w_c})=\bar\Phi\bigl(1/(2\sqrt{w_c^T\Sigma w_c})\bigr)$ we need 
to compare $w_c^T\Sigma w_c$ with 
$w_c^{(1)T}\Sigma\hat w_c^{(1)}$.

To this end let
\begin{align*}
 A_1 & = \{ w : w^T\mu_d=1,\ \| w\|_1\leq c\},
 \\
 A_2 & = \{ w : w^T\hat\mu_d=1,\ \| w\|_1\leq c \}.
\end{align*}
We want to show that for any $w\in A_1$ there exists $\tilde w\in A_2$ 
such that $w^T\Sigma w$ is close to $\tilde w^T\Sigma \tilde w$ 
and vice versa. This means that the minimum 
of $w^T\Sigma w$ over the set 
$A_1$ is close to the minimum over the set $A_2$.

Let $w\in A_1$, and define $\tilde w=w/(w^T\hat\mu_d)$. If 
$\|\tilde w\|_1\leq c$, we have $\tilde w\in A_2$, and  
\[
 w^T\Sigma w = (w^T\hat\mu_d)^2 \tilde w^T\Sigma\tilde w
 =(1+O(c\|\hat\mu_d-\mu_d\|_\infty))^2\tilde w^T\Sigma\tilde w.
\]
If instead $\|\tilde w\|_1> c$, we first define $\bar w\in A_1$ 
and then $w^*=\bar w/(\bar w^T\hat\mu_d)\in A_2$. 
To define $\bar w$ assume without loss of generality that 
$\mu_{d1}=\max_j\{|\mu_{dj}|$. 
Write $w=(w_1,w_{(2)})$ where $w_{(2)}$ is 
$(p-1)$-dimensional, and define $\bar w=(\bar w_1,rw_{(2)})$ 
with $0<r<1$, and $\bar w_1$ chosen such that $\bar w^T\mu_d=1$. 
The latter requirement implies
\[
 \bar w_1\mu_{d1}=1-rw_{(2)}^T\mu_{d(2)}=1-r(1-w_1\mu_{d1}).
\]
We will show that with 
$r=1-c^2\|\hat\mu_d-\mu_d\|_\infty/(c-1/\mu_{d1})
=1-O(c^2\|\hat\mu_d-\mu_d\|_\infty)$ we have 
$\|w^*\|_1\leq c$. 
From the definition of $\bar w$ we have 
\[
 \|\bar w\|_1=|\bar w_1|+r\|\bar w_{(2)}\|_1=
 \frac{|1-r(1-w_1\mu_{d1})|}{\mu_{d1}}+r(\|w\|_1-|w_1|).
\]
If $1-r(1-w_1\mu_{d1})>0$ we get 
\[
 \|\bar w\|_1=\frac{1}{\mu_{d1}}+
 r\bigl(\|w\|_1-\frac{1}{\mu_{d1}}+w_1-|w_1|\bigr)
 \leq \frac{1}{\mu_{d1}}+r\bigl(c-\frac{1}{\mu_{d1}}\bigr).
\]
This shows that $\bar w\in A_1$ and $w^*\in A_2$ since 
\[
 \|w^*\|_1=\frac{\|\bar w\|_1}{\bar w^T\hat\mu_d}\leq 
 \frac{\frac{1}{\mu_{d1}}+r\bigl(c-\frac{1}{\mu_{d1}}\bigr)}{
 1-c\|\hat\mu_d-\mu_d\|_\infty}\leq c,
\]
when $r\leq 1-c^2\|\hat\mu_d-\mu_d\|_\infty/(c-1/\mu_{d1})$. 
If instead $1-r(1-w_1\mu_{d1})<0$ we find 
\[
 \|\bar w\|_1=\frac{r-1}{\mu_{d1}}+ r\|w\|_1 
 \leq rc-\frac{1-r}{\mu_{d1}}\leq rc,
\]
and $\|w^*\|_1\leq rc/(1-c\|\hat\mu_d-\mu_d\|_\infty)\leq c$ for 
$r\leq 1-c\|\hat\mu_d-\mu_d\|_\infty$. The latter condition is 
satisfied with $r\leq 1-c^2\|\hat\mu_d-\mu_d\|_\infty/(c-1/\mu_{d1})$. 
Comparing $\bar w$ and $w$ we get
\begin{align*}
 | w^T\Sigma w-\bar w^T\Sigma\bar w|
 &  \leq 2c\|w-\bar w\|_1 \|\Sigma\|_\infty
 \leq 2c [(1-r)\|w\|_1+(1-r)\frac{1}{\mu_{d1}}] \|\Sigma\|_\infty
 \\
 & =O(c^4\|\hat\mu_d-\mu_d\|_\infty\|\Sigma\|_\infty),
\end{align*}
and also
\[
 | \bar w^T\Sigma\bar w- w^{*T}\Sigma w^*|
  \leq (w^{*T}\Sigma w^*) O(c\|\hat\mu_d-\mu_d\|_\infty).
\]
We have now shown that any value of $w^T\Sigma w$ for $w\in A_1$ is close to 
the corresponding value for some $\tilde w \in A_2$. The other way around, 
starting with $w\in A_2$, is treated in the same way. The only difference is that instead of using $c-1/\mu_{d1}>\epsilon$, we use that when $|\hat\mu_{d1}-\mu_{d1}|<\min\{\epsilon,\epsilon^3/(2+\epsilon^2)\}$, which happens with probability tending to 1 (exponentially fast),
we have $\hat\mu_{d1}>0$ and $c-1/\hat\mu_{d1}>\epsilon/2$. Therefore, the minimum 
$w_c^T\Sigma w_c$ of $w^T\Sigma w$ over the set $A_1$ is close to 
the minimum $w_c^{(1)T}\Sigma w_c^{(1)}$ over the set $A_2$:
\begin{align*}
 w_c^T\Sigma w_c
 & =w_c^{(1)T}\Sigma w_c^{(1)}+
 O(c^4\|\hat\mu_d-\mu_d\|_\infty\|\Sigma\|_\infty)
 +O(c\|\hat\mu_d-\mu_d\|_\infty w_c^T\Sigma w_c)
 \\
 & = w_c^{(1)T}\Sigma w_c^{(1)}+
 O(c^4\|\hat\mu_d-\mu_d\|_\infty\|\Sigma\|_\infty).
\end{align*}
Combining the latter with \eqref{6} we conclude
\[
 |W(\hat \delta)-W(\delta_{w_c})|= 
 O\bigl(c^2 a_n(1+c^2\|\Sigma\|_\infty)\bigr).
\]

\section*{Acknowledgement}

We thank Xin Tong for reading this note and refer to the arXiv version of 
\cite{Fan} for updated versions. 

\appendix

\section*{Appendix A}

An essential step in the proof in \cite{Fan} is the 
inequality (used in equation (21) of that paper)
\[
 \frac{w_c^T\hat\mu_d}{\sqrt{w_c^T\Sigma w_c}}
 \leq 
 \frac{1}{\sqrt{w_c^{(1)T}\Sigma w_c^{(1)}}}.
\]
Unfortunately, this inequality is not correct. We illustrate this by 
a concrete example. We consider the two-dimensional case with 
\[
 \mu_d=(1,0)^T,\ \ 
 \Sigma=\Bigl(\begin{array}{cc} 1 & 1 \\ 1 & \sigma 
 \end{array} \Bigr),\ \ 
 c=1+\epsilon\ \text{with}\ \epsilon<1/\sigma.
\]
In this case we have 
\[
 w_c=(1,-\epsilon)^T\ \text{ and }\ \ w_c^T\Sigma w_c=1-2\epsilon+\sigma\epsilon^2.
\]
Consider next $\hat\mu_d=(1+a,b)$ with $a$ and $b$ small. 
For $a$ and $b$ sufficiently small we obtain
\begin{equation}
 w_c^{(1)}=\bigl(
 \frac{ 1+b[a+\epsilon(1+a)]/(1+a+b) }{1+a},
 -\frac{a+\epsilon(1+a)}{1+a+b}
 \bigr)^T,
 \label{detail1}
\end{equation}
and
\[
 w_c^{(1)T}\Sigma w_c^{(1)}=
 (w_{c1}^{(1)})^2 +2w_{c1}^{(1)}w_{c2}^{(1)}+\sigma (w_{c2}^{(1)})^2.
\]
For $a$ and $b$ small and including $O(a)$ and $O(b)$ terms only 
we get
\begin{equation}
 \frac{1}{\sqrt{w_c^{(1)T}\Sigma w_c^{(1)}}}
 =
 \frac{1}{\sqrt{1-2\epsilon+\sigma\epsilon^2}}
 \bigl\{
 1+a-b\epsilon +(a-b\epsilon)
 \frac{1+\epsilon-\epsilon\sigma(1+\epsilon)}{1-2\epsilon+\sigma\epsilon^2}
 \bigr\},
 \label{11}
\end{equation}
which must be compared to 
\begin{equation}
 \frac{w_c^T\hat\mu_d}{\sqrt{w_c^T\Sigma w_c}}
 =\frac{ 1+a-b\epsilon }{ \sqrt{1-2\epsilon+\sigma\epsilon^2} }.
 \label{12}
\end{equation}
We thus see that \eqref{11} is less that \eqref{12} when $a-b\epsilon$ 
has the opposite sign of $1+\epsilon-\epsilon\sigma(1+\epsilon)$. 
Since $(a-b\epsilon)\sim N(0,(1-2\epsilon+\sigma\epsilon)c_0)$ for 
some constant $c_0$, the probability of a particular sign of 
$a-b\epsilon$ is one half.



\bibliographystyle{ams}

\end{document}